\documentclass[preprint,12pt]{my-elsarticle}

\usepackage[english]{babel}
\usepackage{amsmath}
\usepackage{amssymb}
\usepackage{graphicx}
\usepackage{caption}
\usepackage{url}
\usepackage{a4wide}

\begin{document}

\begin{frontmatter}

\title{Transport and Optimal Control of Vaccination Dynamics for COVID-19}

\author[add1,add2]{Mohamed Abdelaziz Zaitri}
\ead{zaitri@ua.pt}
\ead[url]{https://orcid.org/0000-0002-8956-2266}

\author[add1]{Mohand Ouamer Bibi}
\ead{mobibi.dz@gmail.com}
\ead[url]{https://orcid.org/0000-0002-3531-0928}

\author[add2]{Delfim F. M. Torres\corref{cor1}}
\ead{delfim@ua.pt}
\ead[url]{https://orcid.org/0000-0001-8641-2505}

\address[add1]{Research Unit LaMOS (Modeling and Optimization of Systems),
Department of Operational Research, University of Bejaia, 06000, Bejaia, Algeria}

\address[add2]{Center for Research and Development in Mathematics and Applications (CIDMA),\\
Department of Mathematics, University of Aveiro, 3810-193 Aveiro, Portugal}

\cortext[cor1]{Corresponding author: Delfim F. M. Torres (delfim@ua.pt)}


\begin{abstract}
We develop a mathematical model for transferring the vaccine BNT162b2
based on the heat diffusion equation.
Then, we apply optimal control theory to the proposed generalized SEIR model.
We introduce vaccination for the susceptible population to control the
spread of the COVID-19 epidemic. For this, we use the Pontryagin minimum
principle to find the necessary optimality conditions for the optimal
control. The optimal control problem and the heat diffusion equation
are solved numerically. Finally, several simulations are done to
study and predict the spread of the COVID-19 epidemic in Italy. In particular,
we compare the model in the presence and absence of vaccination.
\end{abstract}

\begin{keyword}
mathematical modeling  \sep COVID-19 pandemic \sep vaccination
\sep optimal control \sep heat diffusion equation.

\MSC[2020] 35K05 \sep 49K15 \sep 92-10.
\end{keyword}

\end{frontmatter}
	

\section{Introduction}

BNT162b2 is an mRNA-based vaccine candidate against SARS-CoV-2, 
currently being developed by Pfizer and BioTech \cite{BNT162b2}. 
As announced on $9^{th}$ November 2020, BNT162b2 shows an 
efficacy against COVID-19 in patients without prior evidence 
of SARS-CoV-2 infection. A first interim efficacy 
analysis was conducted by an external, independent Data Monitoring Committee 
from the Phase 3 clinical study, and the case split, between vaccinated 
individuals and those who received the placebo, indicates a vaccine efficacy 
rate above $90\%$, at seven days after the second dose, of the 94 cases reviewed
\cite{Polack:2020}.

The major obstacle that must be overcome is related to the process 
of transporting the vaccine, which must be stored at $-70^oC$ \cite{Vogel}. 
Pfizer indicates that the vaccine will be distributed from its factories 
in the USA, Belgium and Germany. The American Wall Street Journal revealed 
that Pfizer has developed a special box packed with dry ice and a GPS tracker, 
which can hold 5000 doses of the vaccine under the right conditions for 10 days. 
Moreover, there is another obstacle related to the cost of the transportation boxes, 
where a similar box of 1200 doses in $- 8^oC$ costs 6868 USD, which is very expensive.

The transport of the vaccine must comply with the general standards for drug storage 
and the recommended conditions. Although many transport vehicles are equipped with 
refrigeration devices, assuring recommended storage conditions, simple insulated 
transport boxes are often used. In this study, we use the heat diffusion equation 
and assume that the shape of the vaccine bottle is cylindrical \cite{necati}. 
We perform the calculations to find out an initial temperature that ensures 
the arrival of the vaccine while fulfilling the required condition of $-70^{o}C$, 
by using insulated transfer boxes with the internal temperature at $0^oC$ \cite{Shashkov}.

Optimal control is a mathematical theory that consists of finding a control 
that optimizes a functional on a domain described by a system of differential equations. 
This theory is applied in various fields of the engineering sciences: aeronautics, 
physics, biomedicine, etc. The Pontryagin minimum principle is used to find 
the necessary conditions for optimal controls \cite{Pontr}. 

Several models were presented to predict the spread of COVID-19 
\cite{LemosP,intro,MR4200529,MyID:459,Zamir,Zine}. 
These studies used the SIR and SEIR models \cite{Kermack} 
and the generalized SEIR model \cite{Pengetal}. Most of them 
were implemented to evaluate the strategy of the preventive measures 
\cite{india,Ballesteros,Berhe,aalpha,bbeta,Afrique}.

In \cite{ebola}, the authors present a mathematical model to analyze the Ebola 
epidemic and two optimal control problems related to the transmission 
of Ebola disease with vaccination. In \cite{Pengetal}, the authors present 
a mathematical model to analyze the COVID-19 epidemic based on a dynamic mechanism 
that incorporates the intrinsic impact of hidden latent and infectious cases 
on the entire process of the virus transmission. The authors validate this model 
by analyzing data correlation and forecasting available general data. 
Their model reveals the key parameters of the COVID-19 epidemic. Here, 
we modify the model analyzed in \cite{Pengetal} and consider an 
optimal control problem. More precisely, we introduce an extra variable 
for the number of vaccines used. Secondly, we study the associated optimal 
control problem, solving it numerically. Moreover, in order to find out 
the main parameters, we have performed a numerical simulation of the spread 
of COVID-19 in Italy from $01^{st}$ November 2020 to $31^{th}$ January 2021. 
Finally, we have presented another simulation to find the optimal control, 
and we have compared the models with and without vaccination.

The paper is organized as follows. We begin 
by formulating the vaccination transport model
in Section~\ref{sec2}. In Section~\ref{sec3}, 
we recall the generalized SEIR model. Then, in Section~\ref{sec4}, 
we formulate the generalized SEIR model with vaccination 
as an optimal control problem. The obtained optimal control 
problem is solved numerically in Section~\ref{sec5}.
In Section~\ref{sec6}, we present a discussion concerning 
the spread of COVID-19 in Italy during three months, 
starting from $1^{st}$ November 2020. We end with Section~\ref{sec7} 
of conclusion, including some future research directions.


\section{Vaccine transport model}
\label{sec2}

In this section, we present a model to maintain the effectiveness 
of the vaccine while transporting it from the factory storage area 
to the desired destination. The aim is to know the initial temperature 
that maintains the effectiveness of the vaccine, less than $-70^o$, 
and this by using the available mobile boxes at $0^oC$. Thus, 
we propose the following mathematical model:
\begin{equation} 
\label{heat} 
\left\{
\begin{array}{l l}
\displaystyle
\frac{\partial T(t,x,y,z)}{\partial t } - \alpha \nabla^2T(t,x,y,z)=0,
&  \text{on } \  [0,t_*]  \times\Omega,    \\
T(t_*,x,y,z)=-70^oC,&   \forall (x,y,z)\in \Omega, \\
T(t,x,y,z)=0^oC,&   \forall (t,x,y,z)\in [0,t_*]  \times\partial\Omega,   \\
\Omega=\left\{ (x,y,z)\in \mathbb{R}^3:\ x^2+y^2<r^2,\ 0<z<h   \right\} ,\\
\partial\Omega=\left\{ (x,y,z)\in \mathbb{R}^3:\ x^2+y^2=r^2,\ 0\leq z \leq h \right\} , \\
\end{array}\right.
\end{equation}
where $T(t,x,y,z)$ represents the temperature of the vaccine at the point 
$(x, y, z)$ and the time $t$; $t_*$ is the arrival time of the vaccine; 
and $0^{o}C$ is the temperature inside the box. The sets $\Omega$ 
and $\partial\Omega$ represent the interior and the border 
of the bottle containing the vaccine, respectively, $r$ and $h$ 
are the radius and height of the bottle, respectively, 
and $\alpha$ is the thermal diffusivity defined by
\begin{equation}   
\alpha=\frac{k}{\rho c_{\rho}},
\end{equation}
where $k$ is the thermal conductivity, $c_{\rho}$ is the specific 
heat capacity, and $\rho$ is the density.


\section{Initial mathematical model for COVID-19}
\label{sec3}

The generalized SEIR model proposed by Peng et al. \cite{Pengetal} is
expressed by a seven-dimensional dynamical system defined by
\begin{equation} 
\label{q1} 
\left\{
\begin{array}{l l}
\dot{S}(t) = - \dfrac{\beta S(t) I(t)}{N}-\omega S(t), \\
\dot{E}(t) = \dfrac{\beta S(t) I(t)}{N} - \gamma E(t), \\
\dot{I}(t) =  \gamma E(t)-\delta I(t), \\
\dot{Q}(t) =  \delta I(t)-\lambda(t) Q(t)-\kappa(t) Q(t), \\
\dot{R}(t) = \lambda(t) Q(t), \\
\dot{D}(t) = \kappa(t) Q(t), \\
\dot{P}(t) = \omega S(t), \\
\end{array}\right.
\end{equation}
where the state variables are subjected to the following initial conditions: 
$$
S(0) = S_{0}, \ E(0)=E_0,\ I(0) = I_{0},
\ Q(0) = Q_{0},\ R(0)=R_0, \ D(0) = D_{0}, 
\ P(0) = P_{0}.
$$ 
In this model, the population is divided into the following compartments:  
susceptible individuals $S(t)$, exposed individuals $E(t)$, infected individuals $I(t)$, 
quarantined individuals $Q(t)$,  recovered individuals $R(t)$,  
death individuals $D(t)$, and insusceptible/protected individuals $P(t)$. 
These variables, in total, constitute the whole population, denoted by $N$:  
$$
N = S(t)+E(t)+I(t)+Q(t)+R(t)+D(t)+P(t).
$$
The parameters  $\omega$, $\beta$, $\gamma$, $\delta$, $\lambda(t)$, 
and $\kappa(t)$ represent the protection rate, infection rate, 
inverse of the average latent time, rate at which infectious 
people enter in quarantine, time-dependent recovery rate, 
and the time-dependent mortality rate, respectively. The recovery 
$\lambda(t)$ and mortality $\kappa(t)$ rates  
are analytical functions of time, defined by
\begin{equation} 
\lambda(t)=\frac{\lambda_1}{1+\exp(-\lambda_2(t-\lambda_3))},
\end{equation}
\begin{equation} 
\kappa(t)=\frac{\kappa_1}{\exp(\kappa_2(t-\kappa_3))+\exp(-\kappa_2(t-\kappa_3))},
\end{equation}
where the parameters $\lambda_1$, $\lambda_2$, $\lambda_3$, $\kappa_1$,
$\kappa_2$ and $\kappa_3$ are empirically determined in Section~\ref{sec6}.


\section{Mathematical model for COVID-19 with vaccination}
\label{sec4}

\begin{figure}
\centering
\includegraphics[scale=0.18]{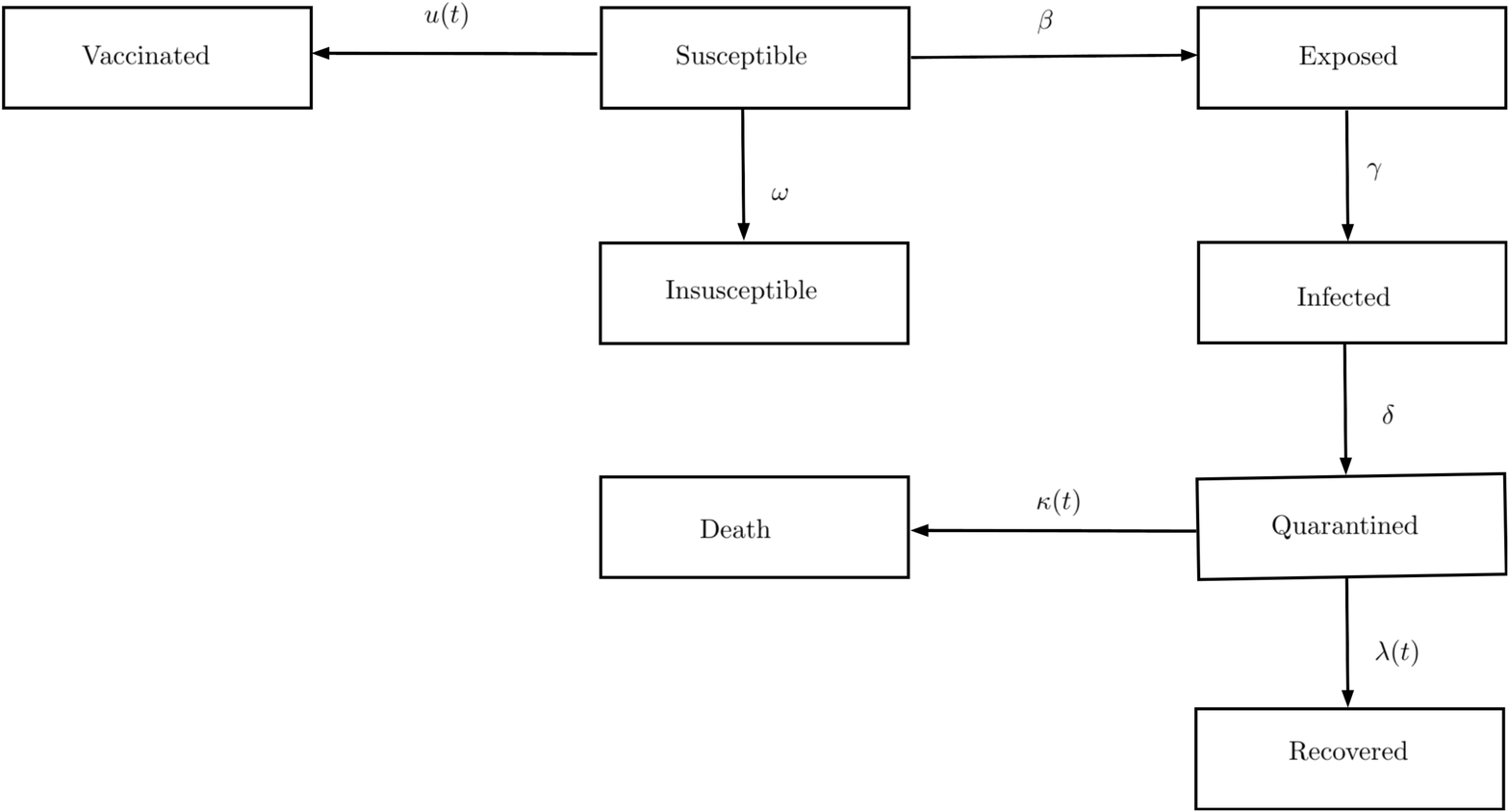}
\caption{Schematic diagram of the generalized SEIR model with vaccination.}
\label{fig10}
\end{figure}

We now introduce the vaccine for the susceptible population in order 
to control the spread of COVID-19. Let us introduce in model \eqref{q1}
a control function $u(t)$ and an extra variable $W(t)$, $t\in [0, t_f ]$, 
representing the percentage of susceptible individuals being vaccinated 
and the number of vaccines used, respectively, with 
\begin{equation} 
\frac{dW}{dt}(t)= u(t)S(t),\ \   
\text{subject to the initial condition} \ W(0)=0,\\
\end{equation}
where $t_f$ represents the final time of the vaccination program. 
Hence, our model with vaccination is given by the following system 
of eight nonlinear ordinary differential equations:
\begin{equation} 
\label{modelv} 
\left\{
\begin{array}{l l}
\dot{S}(t) = - \dfrac{\beta S(t) I(t)}{N}-\omega S(t)-  u(t) S(t), \\
\dot{E}(t) = \dfrac{\beta S(t) I(t)}{N} - \gamma E(t), \\
\dot{I}(t) =  \gamma E(t)-\delta I(t), \\
\dot{Q}(t) =  \delta I(t)-\lambda(t) Q(t)-\kappa(t) Q(t), \\
\dot{R}(t) = \lambda(t) Q(t), \\
\dot{D}(t) = \kappa(t) Q(t), \\
\dot{P}(t) = \omega S(t), \\
\dot{W}(t) =  u(t) S(t),
\end{array}\right.
\end{equation}
where the state variables are subject to the initial conditions: 
\begin{equation*}
\begin{split}
S(0) &= S_{0},\ E(0)=E_0,\ I(0) = I_{0},\ Q(0) = Q_{0},\\
R(0) &=R_0,\ D(0) = D_{0}, \ P(0) = P_{0}, \ W(0)=W_0=0.
\end{split}
\end{equation*}
A schematic diagram of model \eqref{modelv} is given in Figure~\ref{fig10}.

\section{Optimal Control}
\label{sec5}

We consider the model with vaccination (\ref{modelv}) 
and formulate an optimal control problem to determine the vaccination strategy $u$ 
that minimizes the cost of treatment and vaccination:
\begin{equation}
\label{Ju}
\min_{u}J(u)=\int\limits_{0}^{t_f} \left(   w_1 I^2(t)+   w_2 u^2(t) \right) dt,
\end{equation}
where $w_1$ and $w_2$ represent the weights associated with the cost 
of treatment and vaccination, respectively. We assume that the control 
function $u$ takes values between $0$ and $1$. When $ u(t) = 0$, 
no susceptible individual is vaccinated at time $t$ and if $ u(t) = 1$, then
all susceptible individuals are vaccinated at time $t$. 

Let $x(t)=(x_1(t),\ldots,x_8(t))=(S(t),E(t),I(t),Q(t),R(t),D(t),P(t),W(t))\in\mathbb{R}^8$. 
The optimal control problem consists in finding the control $\tilde{u}$ 
and the associated optimal trajectory $\tilde{x}$, satisfying the control 
system (\ref{modelv}) with the given initial conditions
\begin{equation}
\label{ssi}
x(0)=(S_0,E_0,I_0,Q_0,R_0,D_0,P_0,W_0),
\end{equation}
where the control $\tilde{u}\in \Gamma$, 
\begin{equation}\label{GA}
\Gamma=\{u(\cdot)\in L^\infty([0,t_f],\mathbb{R})\ : \  0\leq u(t)\leq 1,
\ t \in [0, t_f]  \},
\end{equation}
minimizes the objective functional (\ref{Ju}).
With the new variables, problem (\ref{modelv})--(\ref{GA}) becomes
\begin{equation}
\label{sssi} 
\begin{gathered}
\min\limits_{u\in\Gamma}
J(u)=\int\limits_{0}^{t_f} \left(   
w_1 x_3^2(t)+   w_2 u^2(t)\right)  dt,\\
\dot{x}(t)=A(t) x(t)+B(x(t))u(t)+f(x(t)),
\quad x(0)=(S_0,E_0,I_0,Q_0,R_0,D_0,P_0,W_0),
\end{gathered}
\end{equation}
where
$$
A(t)=\left( 
\begin{array}{cccccccc}
-\omega & 0 & 0 &0& 0 & 0 &0&0\\
0 & -\gamma & 0 & 0& 0 & 0 &0&0\\
0 & \gamma & -\delta & 0& 0 & 0 &0&0\\
0 & 0 & \delta & -\lambda(t)-\kappa(t)& 0 & 0 &0&0\\
0 & 0 & 0 & \lambda(t)& 0 & 0 &0&0\\
0 & 0 & 0 & \kappa(t)& 0 & 0 &0&0\\
\omega & 0 & 0 &0& 0 & 0 &0&0\\
0 & 0 & 0 &0& 0 & 0 &0&0\\
\end{array}  
\right),
$$
$$
B(x)=( - x_1,\ 0,\ 0, \ 0, \ 0,\ 0 , \ 0 , \   x_1)^T,
$$
$$
f(x)=\left( -\frac{\beta x_1 x_3}{N},\ 
\frac{\beta x_1 x_3}{N},\ 0, \ 0, \ 0,\ 0 , \ 0 , \ 0 \right)^T.
$$
The existence of the optimal control $\tilde{u}$ and the associated optimal trajectory 
$\tilde{x}$ comes from the convexity of the integrand of the cost functional (\ref{Ju}) 
with respect to the control $u$ and the Lipschitz property of the state system with 
respect to the state vector $x(t)$ (see \cite{existence} for existence results of optimal solutions).
According to the Pontryagin Minimum Principle \cite{Pontr}, 
if $\tilde{u}\in\Gamma$  is optimal for the problem (\ref{sssi}) with 
fixed final time $t_f$, then there exists $\psi\in AC([0,t_f],\mathbb{R}^8)$,  
$\psi (t) = (\psi_1(t), \ldots , \psi_8(t))$, called the adjoint vector, such that
$$
\dot{x}=\frac  {\partial   H}{ \partial  \psi}\ 
\text{ and }  \ \dot{\psi}=-\frac{\partial   H}{ \partial   x},
$$
where the Hamiltonian $H$ is defined by
\begin{equation}
H(t,x,\psi, u )=w_1 x_3^2+w_2u^2+\psi^T   \left( A(t)x+B(x)u+f(x) \right).
\end{equation}
The adjoint functions satisfy
\begin{equation}
\label{eq:adj}
\dot{\psi}=-\frac{\partial   H}{ \partial   x}
=\bar{A}(t,x,u)\psi+\bar{B}(x),
\end{equation}
where 
$$
\bar{A}(t,x,u)
=\left( 
\begin{array}{cccccccc}
\frac{\beta x_3 }{N}+\omega+u&  -\frac{\beta x_3}{N} & 0 &0& 0 & 0 &-\omega & - u\\
0 & \gamma & -\gamma& 0& 0 & 0 &0&0\\
\frac{\beta x_1}{N} &- \frac{\beta x_1 }{N}& \delta & -\delta & 0 & 0 &0&0\\
0 &0& 0 & \lambda(t)+\kappa(t)& -\lambda(t)& -\kappa(t) &0&0\\
0 &0& 0 & 0& 0& 0 &0&0\\
0 &0& 0 & 0& 0& 0 &0&0\\
0 &0& 0 & 0& 0& 0 &0&0\\
0 &0& 0 & 0& 0& 0 &0&0
\end{array} \right),
$$
$$
\bar{B}(x)=\left( 0, \ 0 , \ -2w_1x_3, \ 0, \ 0, \ 0, \ 0,\ 0 \right)^T.
$$
The minimality condition
\begin{equation}
\label{ghi}
H(t,\tilde{x}(t),\tilde{\psi}(t),\tilde{u}(t))
=\min\limits_{u\in [0,1]}H(t,\tilde{x}(t),\tilde{\psi}(t),u)
\end{equation}
holds almost everywhere on $[0, t_f ]$. Moreover, the transversality conditions  
assert that $\tilde{\psi}_i(t_f)=0$, $i=1,\ldots,8$. It follows
from the the Pontryagin  minimum principle that the extremal
control $\tilde{u}^p$ is given by
\begin{equation}
\label{eq:utilde}
\tilde{u}^p(t)=\left\{
\begin{array}{ll}
\tilde{u}(t)& \text{if} \ \ 0< \tilde{u}(t) < 1 ,\\
0 & \text{if} \ \  \tilde{u}(t) \leq 0 ,  \\
1 & \text{if} \ \  \tilde{u}(t) \geq 1,\\
\end{array}
\right.
\end{equation}
where
\begin{equation}
\label{utilde}
\tilde{u}(t)=\frac{  \tilde{x}_1(t)  \left(\tilde{\psi}_1(t)-\tilde{\psi}_8(t)\right)}{2w_2}.
\end{equation}


\section{Numerical Results}
\label{sec6}

The current study aims to find the initial temperature to maintain the effectiveness 
of the vaccine during the transportation process, as well as determining 
an optimal vaccination strategy to limit the spread of  COVID-19 in Italy.
For that we reduce the costs of treatment and vaccination, during the three months 
starting from $1^{st}$ November 2020.  We use the MATLAB R2020b program 
to perform all numerical computations. The initial conditions and real data 
are taken from the public database \emph{Dati COVID-19 Italia}, available from
\url{https://github.com/pcm-dpc/COVID-19}.

\begin{figure}
\centering
\includegraphics[scale=0.44]{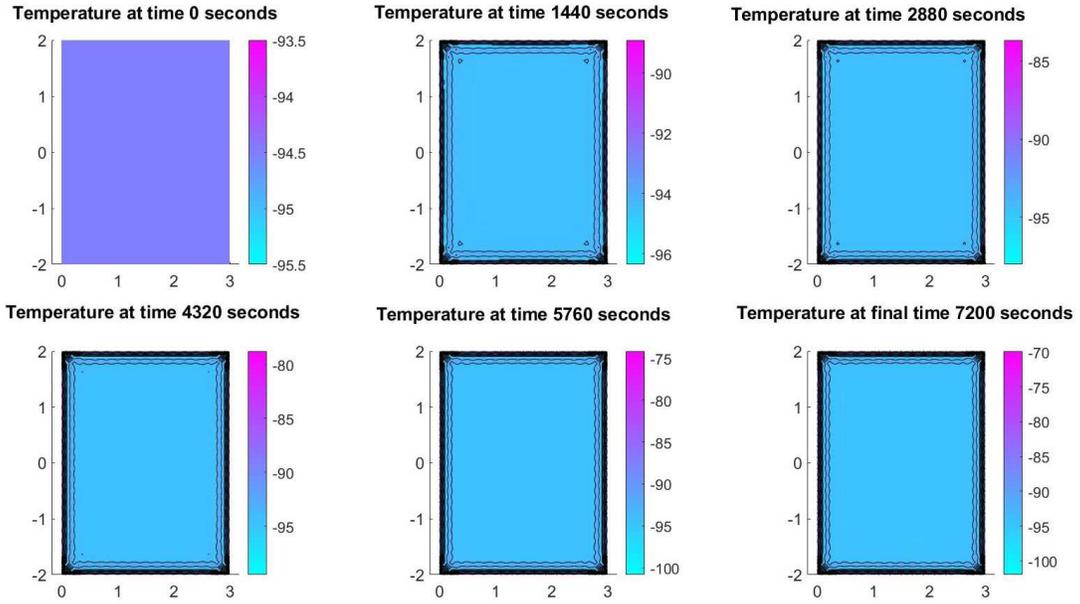}
\caption{Numerical solution of the heat diffusion equation (\ref{heat}).}
\label{Sedp}
\end{figure}

We assume $r=3\ cm$, $h= 4\ cm$, $\alpha = 0.0137 \ W/(m\cdot C)$,
$\rho = 2600 \  kg/m^3$, $c_{\rho }= 750 \  W\cdot s/(kg\cdot C)$ 
and $t_*= 7200\  s$, with the heat transfer coefficient equal to $1$. 
In Fig.~\ref{Sedp} we present the numerical solution of the heat 
diffusion equation \eqref{heat}, which gives the initial temperature equal 
to $-94.5^{o}C$.

We consider the following initial guesses: 
$\omega =0.06$, $\beta=1$, $\gamma=5$, $\delta=0.5$,
$(\lambda_1,\lambda_2,\lambda_3)=(0.01,\ 0.1, \ 10)$ 
and $(\kappa_1,\kappa_2,\kappa_3)=(0.001,\ 0.001,\ 10)$. 
 
The parameters of the generalized SEIR model are computed simultaneously 
by a nonlinear least-squares solver \cite{Cheynet}. These parameters 
over the period starting from $1^{st} $  November 2020 till $31^{th}$ January 2021 are:
$\omega =0.0547$,  $\beta=0.5425$, $\gamma=0.0873$, $\delta=0.3425$, 
$(\lambda_1,\lambda_2,\lambda_3)=(0.0999,\  0.0501,\ 38.8542)$ 
and $(\kappa_1,\kappa_2,\kappa_3)=(0.0021,\ 0.0125,\ 66.6652)$.

\begin{figure}
\centering
\includegraphics[scale=0.28]{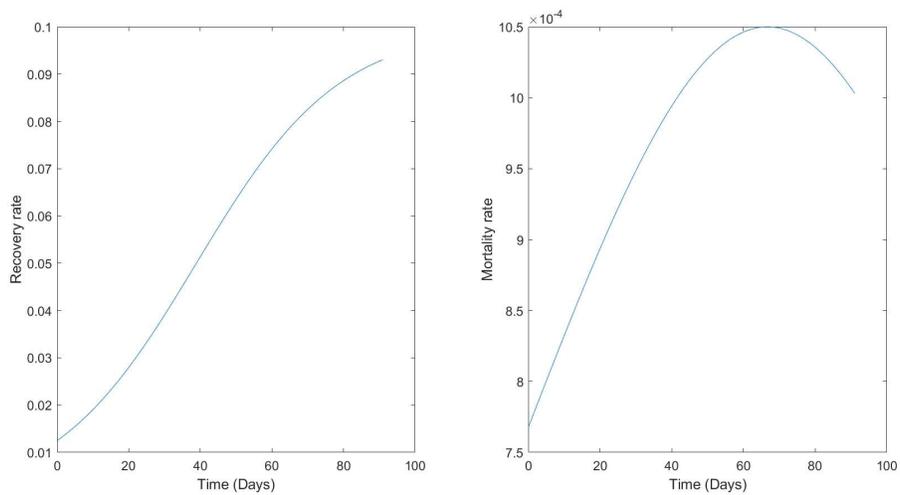}
\caption{The recovery and mortality rates.}
\label{rmm}
\end{figure}

In Fig.~\ref{rmm} we show the recovery rate $\lambda(t)$ and the mortality rate $\kappa(t)$.

\begin{figure}
\centering
\includegraphics[scale=0.38]{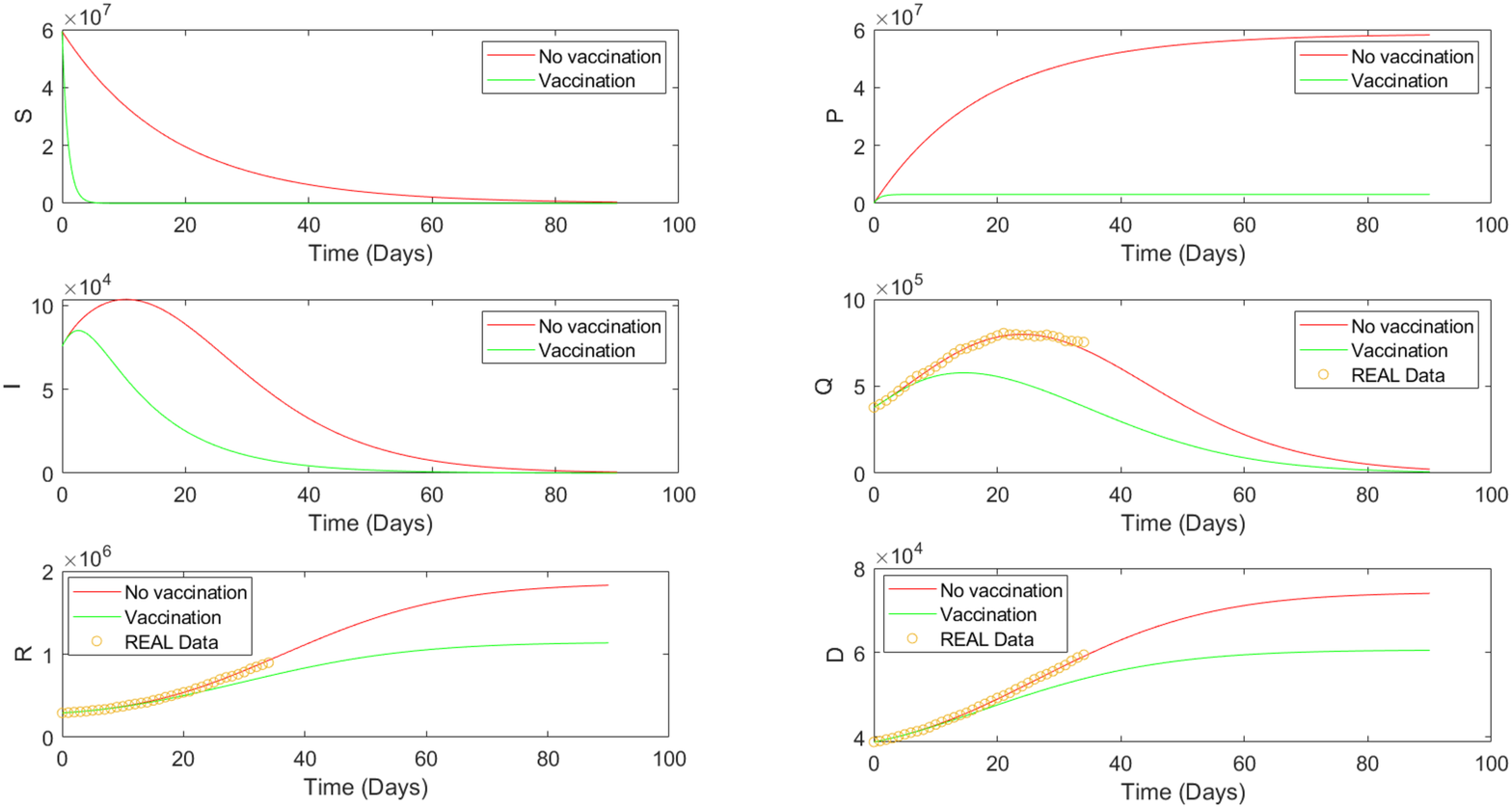}
\caption{The solutions of the generalized SEIR models 
\eqref{q1} and \eqref{modelv}, respectively
without and with vaccination, 
and real data of Italy from $1^{st}$ November till 
$6^{th}$ December 2021 with total population of $N=60.480.000$.}
\label{figs}
\end{figure}

We fixed $w_1=w_2=1$. The numerical solutions to the non-linear differential 
equations that represent the generalized SEIR model \eqref{q1}, 
the generalized SEIR model with vaccination \eqref{modelv},
and the real data of the quarantined, recovered and death cases, 
from $1^{st}$ November till $6^{th}$ December 2021,
are shown in Fig.~\ref{figs}.

\begin{figure}
\centering
\includegraphics[scale=0.27]{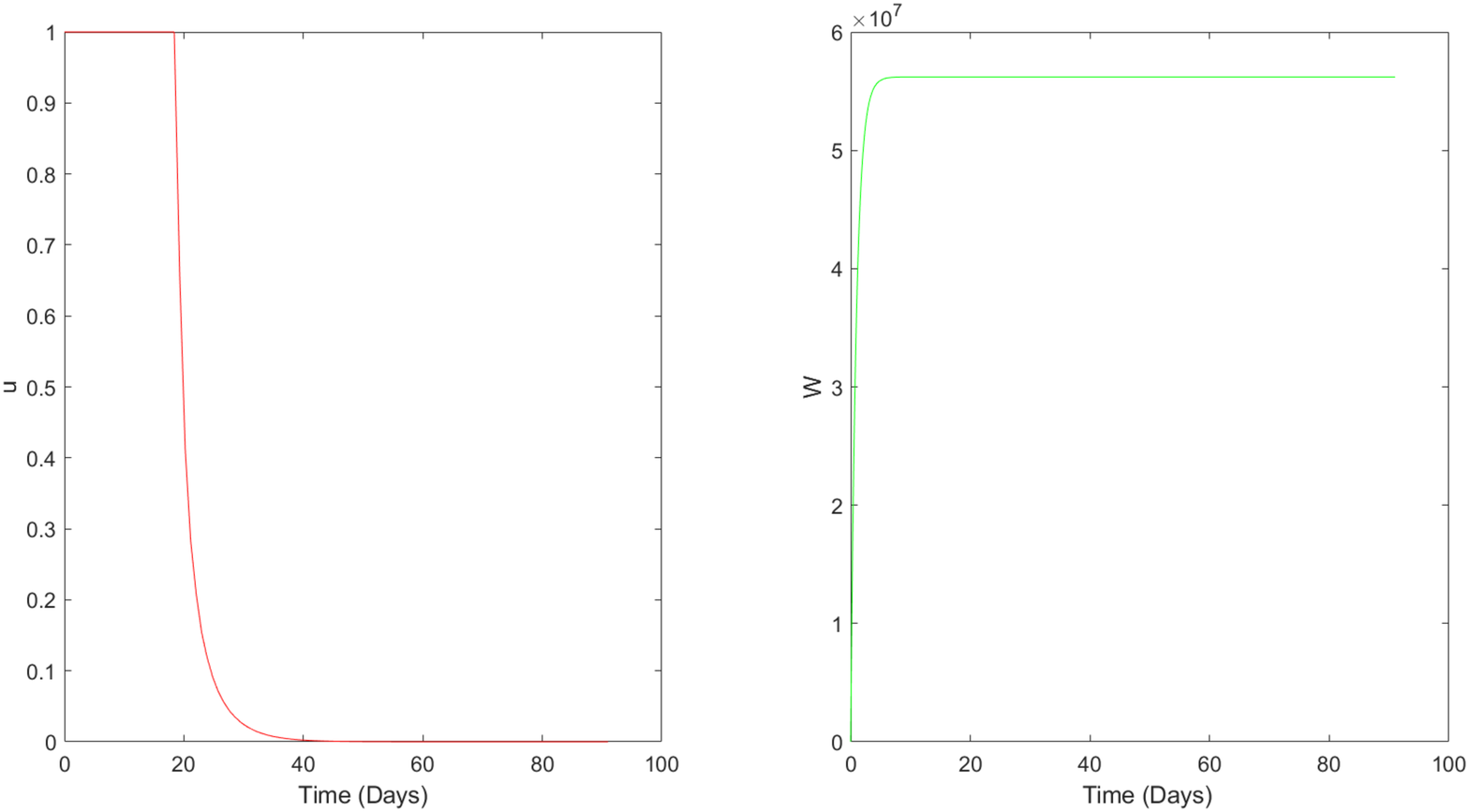}
\caption{The optimal control $\tilde{u}$ (left)
and the number of vaccines $W(t)$ (right).}
\label{figu}
\end{figure}

In Fig.~\ref{figu} we present the optimal control 
\eqref{eq:utilde}--\eqref{utilde}
and the number of vaccines used starting from $ 1^{st} $  
November 2020 till $31^{th}$ January 2021.

The orange curves in Fig.~\ref{figs} represent the real data 
for the number of the quarantine, recovery and death cases 
in Italy starting from $1^{st}$ November till $6^{th}$ December $2020$. 
The red curves in Fig.~\ref{figs} represent the solutions of 
the generalized SEIR model \eqref{q1} without vaccination, 
and they simulate what happen from the beginning of November 
to the end of January. There is an increase in the number of the recovered, 
death and insusceptible cases that reach, respectively, 
$1.830.000$, $74.050$ and $58.130.000$ cases. The red curves 
for both the number of infected and quarantined individuals 
have their higher limit values of $103500$ cases on $11^{th}$ November 
and $798.500$ on $25^{th}$ November, respectively, reaching the values 
$614$ and $22.640$ cases on $31^{th}$ January $2021$, respectively. 
We note that the number of susceptible individuals gradually decrease, 
reaching $416.600$ cases at the end of January $2021$.

The green curves in Fig.~\ref{figs} represent the solutions 
of the generalized SEIR model \eqref{modelv} with vaccination, 
and they simulate what happened from the beginning of November 
to the end of January. There is an increase in the number of 
recovered, death and insusceptible cases that reach, respectively, 
$1.135.000$, $60.560$  and $3.076.000$ cases. The green curves 
for both the number of infected and quarantined individuals 
have their higher limit values of $84.800$ cases on $4^{th}$ November 
and $577.600$ cases on $15^{th}$ November, respectively, reaching 
$55$ and $7.237$ cases on $31^{th}$ January $2021$, respectively. 
We note that the number of susceptible individuals decrease rapidly 
reaching $0$ cases on $19^{th}$ November $2020$.

The red curve in Fig.~\ref{figu} shows that the optimal vaccination 
of 100 percent of the susceptible individuals takes 19 days, 
followed by a rapid decrease in the number of susceptible individuals, 
which means they move to the class of vaccinated. The green curve in Fig.~\ref{figu} 
shows the necessary number of vaccines to eliminate COVID-19, 
which is estimated at 56.200.100 doses. The total number of vaccinated 
and insusceptible individuals equal to 59.276.100 
of the total Italian population of 60.480.000.


\section{Conclusion}
\label{sec7}

Our results show the importance of the vaccine for COVID-19 control
and also the best result that could be obtained if the number 
of available vaccines satisfies the needs of the population
and are distributed according with the theory of optimal control.

Here our optimal control problem has only one control: the vaccine.
In reality, there are several other factors to take into account
and other variables to control. In a future work, we would like 
to use the support maximum principle  \cite{bibi3,bibi2,fouzia}, 
as well as the hybrid direction method \cite{zaitri}, to elaborate 
a primal-dual method for solving a more realistic optimal control problem, 
in presence of multiple inputs \cite{bibi}.


\section*{Acknowledgments}

This research is part of first author's Ph.D. project.
Zaitri is grateful to the financial support from the 
Ministry of Higher Education and Scientific Research of Algeria;
Torres acknowledges the financial support from
CIDMA through project UIDB/04106/2020.



\end{document}